\documentclass[reqno,centertags, 11pt]{amsart}
\usepackage{amscd}
\usepackage{latexsym}
\usepackage{amsmath}
\usepackage{amsthm}
\usepackage{amsfonts}
\usepackage{amssymb}
\numberwithin{equation}{section}

\newcommand{\T}{\partial\mathbb{D}}

\newcommand{\ds}{\displaystyle}
\newcommand{\vp}{\varphi}
\newcommand{\h}{h}
\newcommand{\p}{\psi}
\newcommand{\s}{s}
\newcommand{\kt}{\tilde{K}}
\newcommand{\tT}{\widetilde{T}}
\newcommand{\ttau}{\tilde{\tau}}
\newcommand{\dt}{\tilde{\delta}}
\newcommand{\rt}{\tilde{\rho}}
\newcommand{\ol}{\overline}
\newcommand{\supp}{\text{\rm{supp}}}
\newcommand{\dist}{\text{\rm{dist}}}

\newtheorem{lemma}{Lemma}[section]
\newtheorem{theorem}{Theorem}[section]

\begin{document}
\title[$1^{st}$ AND $2^{nd}$ TYPE PARAORTHOGONAL POLYNOMIALS AND THEIR ZEROS]{First and Second Kind Paraorthogonal Polynomials and their Zeros}
\author[M.-W. L. Wong]{Manwah Lilian Wong*}
\thanks{$^*$ MC 253-37, Mathematics Department, California Institute of Technology, Pasadena, CA 91125, USA.
E-mail: wongmw@caltech.edu. Supported by the Croucher Foundation Scholarship, Hong Kong}

\date{April 17th, 2006} 

\begin{abstract} Given a probability measure $\mu$ with infinite support on the unit circle $\partial\mathbb{D}=\{z:|z|=1\}$, we consider a sequence of paraorthogonal polynomials $\h_n(z,\lambda)$ vanishing at $z=\lambda$ where $\lambda \in \T$ is fixed. We prove that for any fixed $z_0 \not \in \supp(d\mu)$ distinct from $\lambda$, we can find an explicit $\rho>0$ independent of $n$ such that either $\h_n$ or $\h_{n+1}$ (or both) has no zero inside the disk $B(z_0, \rho)$, with the possible exception of $\lambda$. 

Then we introduce paraorthogonal polynomials of the second kind, denoted $\s_n(z,\lambda)$. We prove three results concerning $\s_n$ and $\h_n$. First, we prove that zeros of $\s_n$ and $\h_n$ interlace. Second, for $z_0$ an isolated point in $\supp(d\mu)$, we find an explicit radius $\rt$ such that either $\s_n$ or $\s_{n+1}$ (or both) have no zeros inside $B(z_0,\rt)$. Finally we prove that for such $z_0$ we can find an explicit radius such that either $\h_n$ or $\h_{n+1}$ (or both) has at most one zero inside the ball $B(z_0,\rt)$.

\end{abstract}
\maketitle
\section{Introduction}

Suppose we are given a probability measure $\mu$ on the unit circle $\T=\{z\in \mathbb{C}: |z|=1\}$ with infinite support. We form the inner product $\langle,\rangle$ and the norm in $L^2(d\mu)$ as follows:
\begin{equation} \begin{array}{lllllll}
\langle f,g \rangle & = & \ds \int_{\T}{f(z)}\overline{g(z)} d\mu(z) & ; &  \|f\| & = & \langle f,f \rangle^{1/2}\, .
\end{array} \end{equation}
By the Gram-Schmidt process, we then obtain a sequence of monic orthogonal polynomials $(\Phi_n)_{n=1}^{\infty}$, the normalized sequence being $(\varphi_n)_{n=1}^{\infty}$, such that $\vp_n$ is an $n^{th}$ degree polynomial with the property:
\begin{equation} 
\langle \vp_m, \vp_n \rangle =  \delta_{mn} \, .
\end{equation}
These orthogonal polynomials satisfy the Szeg\H o recursion relation:
\begin{equation}
\Phi_{n}(z) = z\Phi_{n-1}(z)-\overline{\alpha_{n-1}}\Phi_{n-1}^*(z)
\end{equation}
where $\Phi_m^*(z)=z^m \overline{\Phi_m(1/\overline{z})}$.

The family of $\alpha_n$'s are known as the Verblunsky coefficients. There are a few important properties of orthogonal polynomials and Verblunsky coefficients which are relevant to this paper:
\begin{align}
&|\alpha_n|<1 \\
& \|\Phi_n\| = (1-|\alpha_{n-1}|^2)^{1/2}\|\Phi_{n-1}\| = \ds \prod_{j=0}^{n-1}(1-|\alpha_j|^2)^{1/2} \\
& |\Phi_n(z)| =|\Phi_n^*(z)| \Leftrightarrow z\in \T  \\
& \Phi_n(z) \text{ has all its zeros inside  } \mathbb{D}  \\
& \langle \Phi_n(z),z^k \rangle = 0 \; \mbox{for} \; k =0,\dots, n-1 \label{orthogonalities1}\\
& \langle \Phi_n^*(z),z^k \rangle = 0 \; \mbox{for} \; k =1,\dots, n \, .
\label{orthogonalities2}
\end{align}
\medskip

\emph{Paraorthogonal polynomials} were introduced at least as early as in \cite{jones}. An $n^{th}$ degree paraorthogonal polynomial is of the form (up to  multiplication with a constant): \begin{equation}
H_{n}(z, \beta_{n-1})= z\Phi_{n-1}(z)-\ol{\beta_{n-1}} \Phi_{n-1}^*(z)
\label{popdef}
\end{equation}
with $\beta_{n-1} \in \T$; $\Phi_{n-1}^*(z)=z^{n-1} \overline{\Phi_{n-1}(1/{\overline{z}})}$. \\

Paraorthogonal polynomials have a lot in common with orthogonal polynomials on the real line $(p_n)_{n=0}^{\infty}$. For instance, a paraorthogonal polynomial has simple zeros on the unit circle while $p_n$ has simple zeros on the real line. Besides, for a specific family of paraorthogonal polynomials $(h_n)_{n=0}^{\infty}$ that we shall consider, it has been proven in \cite{cmv, golinskii} that zeros of $h_n$ and $h_{n+1}$ strictly interlace, this interlacing property is also shared by $p_n$ and $p_{n+1}$.
\medskip

In this paper we shall prove three results concerning this specific family of paraorthogonal polynomials $(h_n)_{n=0}^{\infty}$, namely Theorem \ref{theorem1}, Theorem \ref{theorem2} and Theorem \ref{theorem3}. These results are in parallel with those proven for orthogonal polynomials of the real line $p_n$.

Theorem \ref{theorem1} and Theorem \ref{theorem3} are analogues of the following results by Denisov--Simon \cite{denisov--simon}:

\begin{theorem}
Let $\delta=\dist(x_0, \supp(d\mu))>0$. Suppose $a_{n+1}$ is the recursion coefficient as given by $xp_n(x)=a_{n+1} p_{n+1}(x)+b_{n+1}p_n(x)+a_n p_{n-1}(x)$. Let $r_n=\delta^2/(\delta+\sqrt{2}a_{n+1})$. Then either $p_n$ or $p_{n+1}$ (or both) has no zeros in $(x_0-r_n,x_0+r_n)$.
\end{theorem}

\begin{theorem} 
Let $x_0$ be an isolated point of $\supp(d\mu)$ on the real line. Then there exists $d_0>0$ so that if $\delta_n=d_0^2/(d_0+\sqrt{2}a_{n+1}))$, then at least one of $p_n$ and $p_{n+1}$ has no zeros or one zero in $(x_0-\delta_n,x_0+\delta_n)$.
\end{theorem}
\medskip

Theorem \ref{theorem2}, which proves that first and second kind paraorthogonal polynomials of the same degree have interlacing zeros, is an analogue of the following well-known fact about first and second kind orthogonal polynomials on the real line, $p_n$ and $q_{n}$:
\begin{theorem} Zeros of $p_n$ and $q_n$ strictly interlace.
\end{theorem}

For a more comprehensive introduction to orthogonal polynomials and paraorthogonal polynomials, the reader should refer to \cite{onefoot, simon, szego}. \\

\begin{section}{properties of paraorthogonal polynomials}\label{propsec}

A major difference between orthogonal polynomials and paraorthogonal polynomials lies in the fact that $\alpha_n \in \mathbb{D}$ is determined uniquely by the measure, while $\beta_n \in \T$ could be chosen arbitrarily on the unit circle. These differences give rise to the following properties of $H_n$ which are not shared by $\Phi_n$:\\

\emph{1. Zeros on $\T$ } \quad Unlike orthogonal polynomials which have zeros strictly inside the unit disk, paraorthogonal polynomials have zeros in $\T$. To see that it suffices to note that
\begin{equation}
\ds \left| \frac{z \Phi_n(z)}{\Phi_n^*(z)} \right|=1 \Leftrightarrow z \in \T \, .
\end{equation}
\medskip

\emph{2. Orthogonality} \quad An $n^{th}$ degree paraorthogonal polynomial is orthogonal to $\{z,z^2, \dots, z^{n-1}\}$ because of the orthogonal properties of $\Phi_{n-1}$ and $\Phi_{n-1}^*$ as in (\ref{orthogonalities1}) and (\ref{orthogonalities2}). However, we note that $H_n$ is never orthogonal to $1$ or $z^n$ because
\begin{align}
\langle H_n,1\rangle&=(\ol{\alpha_{n-1}}-\ol{\beta_{n-1}})\|\Phi_{n-1}\|^2\not = 0 \, ,\\
\langle H_n,z^n\rangle&=\left( 1-\ol{\beta_{n-1}}\alpha_{n-1} \right)\|\Phi_{n-1}\|^2 \not = 0
 \, .
\end{align} \medskip

\emph{3. Representation} \quad Suppose $\lambda$ is a zero of $H_{n}(z, \beta_{n-1})$. We prove that $H_{n}$ could be represented using the reproducing kernel $K_{n}(z,\lambda)= \sum_{j=0}^{n}\vp_j(z)\overline{\vp_j(\lambda)}$ and a constant C as follows: 
\begin{equation} H_{n}(z,\beta_{n-1})=\ds 
C(z-\lambda)\ds \sum_{j=0}^{n-1}\varphi_j(z) \overline{\varphi_j(\lambda)}=C(z-\lambda)K_{n-1}(z,\lambda)  \, .
\label{formulah} 
\end{equation}

The argument is related to Szeg\H o \cite{szego} when he proved the Christoffel--Darboux formula. It goes as follows: since $\lambda$ is a zero of $H_{n}$,  $H_{n}(z)=(z-\lambda)h(z)$ for some polynomial $h$ of degree $n-1$. By the orthogonality of $H_{n}$ against $\{z,\dots,z^{n-1}\}$, $\langle zh, z^m \rangle= \langle \lambda h,z^m \rangle$ for $1 \leq m \leq n-1$, which implies that $\ol{\lambda} \langle h,z^{m-1} \rangle= \langle h, z^m \rangle$. Applying this formula recursively, we conclude that 
\begin{equation}
\langle h,z^m\rangle=\ol{\lambda^m}\langle h,1\rangle, \text{ for }0 < m \leq n-1 \, .
\end{equation}
When $m=0$ the argument is trivial. If $\vp_s(z)= \sum_{j=0}^{s}a_j z^j$, then for $0 \leq s\leq n-1$,
\begin{equation}
\langle h,\vp_s\rangle = \langle h,1\rangle \ds\sum_{j=0}^{s} \ol{a_j\lambda^j}  = \langle h,1\rangle \ol{\vp_s(\lambda)}   \, .
\end{equation}

If we express $h$ using Fourier series,
\begin{equation}
h(z) = \ds \sum_{j=0}^{n-1} \langle h,\vp_j\rangle  \vp_j(z) = \ds \sum_{j=0}^{n-1} \langle h,1\rangle  \ol{\vp_j(\lambda)} \vp_j(z) = \langle h,1\rangle K_{n-1}(z,\lambda)\label{hn} \, .
\end{equation}

\emph{4. Simple Zeros}
Let $\lambda$ and $h$ be defined as above. By (\ref{hn}), $\langle h, 1\rangle =0$ implies $h=0$, hence $\langle h, 1\rangle \not = 0$ . In addition, $\vp_0=1$ implies $K_{n-1}(\lambda,\lambda)>0$. Therefore $h(\lambda)=\langle h,1 \rangle K_{n-1}(\lambda,\lambda) \not = 0$. This shows that zeros of paraorthogonal polynomials are simple. \\
\medskip

\emph{5. Linear Independence} \quad The argument for property (3) above also tells us that a paraorthogonal polynomial could vanish at one arbitrary point on the unit circle, and that particular zero fixes the remaining ones. Therefore, two paraorthogonal polynomials of the same degree are linearly independent if and only if all their zeros are distinct.\\
\medskip

The reader could refer to \cite{cmv, simon} for more properties of paraorthogonal polynomials. \\
\medskip
\end{section}

\section{Equivalent Definitions of $h_n$}

Fix $\lambda \in \T$. We define the family of paraorthogonal polynomials $(h_n(z,\lambda))_n$ as follows:
\begin{equation}
\h_n(z, \lambda):= (1-\ol{\lambda}z) K_{n-1}(z,\lambda) \, .
\end{equation}

We will soon see that there are three equivalent definitions of $h_n$ by the Christoffel-Darboux formula. The formula says that for $ \ol{y}z\not = 1$, the reproducing kernel $K_{n-1}(z,y)$ could be expressed in the following ways:
\begin{eqnarray}
K_{n-1}(z,y) & = & \ds \frac{\overline{\vp_{n}^*(y)}\vp_{n}^*(z)-\overline{\vp_{n}(y)}\vp_{n}(z)}{1-\overline{y}z} \\ & = & \ds \frac{\overline{\vp_{n-1}^*(y)}\vp_{n-1}^*(z)-\overline{y}z\overline{\vp_{n-1}(y)}\vp_{n-1}(z)}{1-\overline{y}z} \label{zlambdaexp} \, .
\end{eqnarray}

Hence, we have the following three equivalent definitions of $h_n(z,\lambda)$:
\begin{eqnarray}
\h_n(z) & = & (1-\ol{\lambda} z)\ds \sum_{j=0}^{n-1} \vp_j(z)\ol{\vp_j(\lambda)}
 \label{hndef3} \\
& = & \overline{\vp_{n}^*(\lambda)}\vp_{n}^*(z)-\overline{\vp_{n}(\lambda)}\vp_{n}(z) \label{hndef1}\\
& = & \overline{\vp_{n-1}^*(\lambda)}\vp_{n-1}^*(z)- z\overline{\lambda} \overline{\vp_{n-1}(\lambda)}\vp_{n-1}(z) \label{hndef2} \, . 
\end{eqnarray}

By rewriting (\ref{hndef2}) in the form of (\ref{popdef}),
\begin{equation}
\h_n(z) = -\ol{\lambda \vp_{n-1}(\lambda)}\left( z\vp_{n-1}(z)-\ds \lambda\frac{ \ol{\vp_{n-1}^*(\lambda)}}{\ol{\vp_{n-1}(\lambda)}} \vp_{n-1}^*(z) \right)
\label{hnspecial}
\end{equation}
we see that the coefficient $\beta_{n-1}$ of this particular family of paraorthogonal polynomials are
\begin{equation}
\beta_{n-1}(h_n)=\ds \ol{\lambda}\frac{ \vp_{n-1}^*(\lambda)}{\vp_{n-1}(\lambda)} \, .
\end{equation}
\medskip

\section{Paraorthogonal Polynomials of the Second Kind $s_n$}

Paraorthogonal polynomials of the second kind arise from orthogonal polynomials of the second kind, namely $\p_k(z)$, which are orthogonal polynomials associated to the measure $\nu$ with Verblunsky coefficients
\begin{equation}
\alpha_n(d\nu)= - \alpha_n(d\mu) \, .
\end{equation}
The existence of the measure is guaranteed by Verblunsky's theorem which says that for any given sequence of complex numbers inside $\mathbb{D}$, there corresponds a measure on the unit circle with such as Verblunsky coefficients. \\

With the same $\lambda$ as we used to define $\h_n(z,\lambda)$, we define our \emph{Paraorthogonal Polynomials of the Second Kind} $s_n$ as follows:
\begin{equation}
\s_n(z)=\ol{\vp_{n-1}^*(\lambda)}\p_{n-1}^*(z)+z\ol{\lambda}\ol{\vp_{n-1}(\lambda)}\p_{n-1}(z)  \, .
\label{spopdef}
\end{equation}

If we rewrite (\ref{spopdef}) in the form of (\ref{hnspecial})
\begin{equation}
\s_n(z)=\ol{\lambda \vp_{n-1}(\lambda)}\left( z\p_{n-1}(z)+\ds \lambda\frac{ \ol{\vp_{n-1}^*(\lambda)}}{\ol{\vp_{n-1}(\lambda)}} \p_{n-1}^*(z) \right)
\end{equation}
we see that the $\beta_n$ coefficient of this family of paraorthogonal polynomials $(s_n)_n$ is given by: 
\begin{equation}
\beta_n(\s_n)=-\beta_n(\h_n)  \, .
\end{equation}
\medskip

As in the case of $h_n$, we shall see that there are three equivalent definitions of $s_n$ by means of the \emph{Mixed Christoffel--Darboux Formulae}, which state that:
\begin{gather}
\ol{\vp_{n-1}^*(y)}\p_{n-1}^*(z)+z\ol{y}\ol{\vp_{n-1}(y)}\p_{n-1}(z) = \ol{\vp_{n}^*(y)}\p_{n}^*(z)+\ol{\vp_{n}(y)}\p_{n}(z)
\label{mixformula1}  \; \\
\ds \sum_{j=0}^{n-1} \overline{\vp_j(y)}\p_j(z) =  \frac{2-\overline{\vp_{n}^*(y)}\p_{n}^*(z)-\overline{\vp_{n}(y)}\p_{n}(z)}{1-\overline{y}z} \quad \mbox{ for }y \not = z \, .
\label{mixformula2}
\end{gather}
The reader should refer to Chapter 3.2 of \cite{simon} for the proof. \\

By (\ref{mixformula1}) and (\ref{mixformula2}), $s_n(z,\lambda)$ has the following three equivalent definitions: 
\begin{align}
\s_n(z) & = \ol{\vp_{n-1}^*(\lambda)}\p_{n-1}^*(z)+z\ol{\lambda}\ol{\vp_{n-1}(\lambda)}\p_{n-1}(z)\label{sndef2} \\
& = \ol{\vp_{n}^*(\lambda)}\p_{n}^*(z)+\ol{\vp_{n}(\lambda)}\p_{n}(z) \label{sndef1} \\
& = -(1-\ol{\lambda} z)\ds \sum_{j=0}^{n-1} \overline{\vp_j(\lambda)}\p_j(z)+2  \, . \label{sndef3}
\end{align}

\section{Results}

We prove four results concerning $\h_n$, $\h_{n+1}$, $\s_n$ and $\s_{n+1}$. Some related results will be discussed.\medskip

\begin{theorem}\label{theorem1}
Suppose $z_0 \in \T$ distinct from $\lambda$ and $\delta=\dist(z_0,\supp(d\mu))>0$. Then in the open disk around $z_0$ with radius
\begin{equation}\rho=\frac{\delta^3}{8+\delta^2}
\label{rho}
\end{equation}
either $\h_n$ or $\h_{n+1}$ \rm{(}or both\rm{)} has no zero inside, with the possible exception of $\lambda$.

Furthermore, if $L=\dist(\lambda,\supp(d\mu))>0$, then the radius could be taken as:
\begin{equation} \rho'=\ds \frac{\delta^2 L}{8+\delta L}\label{rho1}  \, .
\end{equation} 
\end{theorem}
Note that when $L > \delta$, $\rho' > \rho$, hence (\ref{rho1}) improves (\ref{rho}).
\vspace{0.5cm}

There is a related conjecture concerning double limit points which was proposed in \cite{golinskii} and proven in \cite{cmv1}. The result says that the set of double limit points of $\h_n$ coincides with $\supp(d\mu)$, except at most the point $\lambda$. In other words, if $\dist(z_0,\supp(d\mu))>0$, then for any sequence of integers $I$, there exists a subsequence $I'\subset I$ and $\epsilon_I>0$ such that for $n \in I'$, either $\h_n$ or $\h_{n+1}$ (or both) has no zero in the open disk $B(z_0, \epsilon_I)$.

However, Theorem \ref{theorem1} is clearly stronger because we found an explicit  radius $\rho$ for which the double zero result holds (\ref{rho}) and the result does not depend on $n$.

\medskip

\begin{theorem} The zeros of $\h_n$ and $\s_n$ strictly interlace, that is, between any two zeros of $\h_n$ $($or $\s_n$$)$, there is one and only one zero of $\s_n$ $($or $\h_n$ respectively$)$ in between.
\label{theorem2}
\end{theorem}

At the same time that this result was proven,  Simon \cite{simon1} demonstrated another way of proving the result using the theory of rank one perturbations of unitary operators. He made the observation that the CMV matrix associated to $\s_n$ is just the original one with the signs of $\alpha_j$ and $\beta_{n-1}$ reversed, and it is unitarily equivalent to one where the signs are not reversed but the first column has opposite sign.

The main tools of the proof are the two real-valued functions $\sigma_n$ and $\eta_n$ which we will define in (\ref{sigman}) and (\ref{etan}). They were used in \cite{cmv} to prove that zeros of $h_n$ and $h_{n+1}$ interlace, but the method employed in the our proof is different. \\

The remaining two results are:
\begin{lemma} \label{mainlemma} Suppose $z_0$ is an isolated point in $\supp(d\mu)$. Then 
\begin{equation}\dt=\dist(z_0,\supp(d\nu))>0 \label{dt}
\end{equation}
and in the ball around $z_0$ with radius
\begin{equation}
\rt=\ds \frac{\dt^2 |z_0-\lambda|}{8+|z_0-\lambda|\dt}
\label{rad}
\end{equation}
either $\s_n$ or $\s_{n+1}$ $($or both$)$ has no zeros inside.
\end{lemma}
\smallskip

\begin{theorem} Suppose $z_0$ is an isolated point of $\supp(d\mu)$ and $\dt$ is as defined in (\ref{dt}). Then in the open disk around $z_0$ with radius
\begin{equation}
\rt= \frac{\dt^2 |z_0-\lambda|}{8+|z_0-\lambda|\dt}
\end{equation}
either $\h_m$ or $\h_{n+1}$ $($or both$)$ has at most one zero inside.
\label{theorem3} \end{theorem}
\medskip

\section{Proof of Theorem \ref{theorem1}}

Before we start the proof, we refer to a theorem about zeros of $\h_n$ in a gap of the measure:
\begin{theorem}(Corollary 2 of \cite{cmv}, Theorem 2 of \cite{golinskii}, Theorem 2.3 of \cite{simon1}) Let an arc $\Gamma=(\alpha, \beta)$ on $\T$ be a gap in $\supp(d\mu)$, that is, $\supp(d\mu)\cap \Gamma=\emptyset$ and $\alpha$ goes to $\beta$ counterclockwise. Then for each $n$, the paraorthogonal polynomial $\h_n$ has at most one zero in $\overline{\Gamma}=[\alpha,\beta]$. \smallskip
\label{golinskii}\end{theorem} 

If $\lambda$ is in a gap $\Gamma$, since $\lambda$ is zero of all $\h_n$, by Theorem \ref{golinskii} above there are no other zeros of $\h_n$ or $\h_{n+1}$ in $\Gamma$. In other words, if $z_0$ and $\lambda$ are in the same gap, in a radius $\delta=\dist(z_0,\supp(d\mu))$ around $z_0$ there could be no zeros other than $\lambda$. Since $\delta>\rho$, Theorem \ref{theorem1} holds. Hence if $\lambda$ is in a gap, it suffices to look at the case when $z_0$ that sits in gaps other than $\Gamma$. In such a situation, $|z_0-\lambda| \geq \dist(z_0,\supp(d\mu))$.

However, if $\lambda$ is not in a gap, that is, $\lambda$ is in the support of a measure, then clearly $|z_0-\lambda|\geq \dist(z_0,\supp(d\mu))$.

Without loss of generality, we may assume that $|z_0-\lambda| \geq \delta$ in this section. \\

We shall divide the proof into two lemmas:\\
\begin{lemma}\label{lemma1}
\begin{equation}  \ds \left|\frac{\h_i(z_0)}{K_{n-1}(z_0,z_0)^{1/2}} \right| \geq \frac{1}{4}|\vp_n(\lambda)| \delta^2
\end{equation}
where 
$i= \begin{cases} n & \mbox{ if } |\h_{n+1}(z_0)| \leq |\h_n(z_0)| \\
n+1 & \mbox{ if } |\h_{n}(z_0)| \leq |\h_{n+1}(z_0)|
\end{cases}$\, .

\end{lemma}
\smallskip 
\begin{proof} Suppose $|\h_{n+1}(z_0)| \leq |\h_n(z_0)|$.

First, we give a bound for the $L^2(\mu)$ norm of $\|(z_0- \cdot)K_{n-1}(z_0, \cdot)\|$. \\
By the parallelogram equality and the fact that $|\vp_n^*(z_0)|=|\vp_n(z_0)|$,
\begin{equation} \begin{array} {lll}
& \|(z_0- \cdot)& K_{n-1}(z_0, \cdot)\|^2 \\&   &= \|\overline{\vp^*_n(\cdot)} \vp^*_n(z_0)-\overline{\vp_{n}(\cdot)}\vp_n(z_0)\|^2 \\
& &\leq  2 |\vp_n^*(z_0)|^2+2 |\vp_n(z_0)|^2 \\
&  & =
4 \left| \ds \frac{\h_{n+1}(z_0)-\h_n(z_0)}{(z_0-\lambda)\overline{\vp_n(\lambda)}}\right|^2 \smallskip\\
& &\leq  \ds \frac {4|\h_{n+1}(z_0)|^2+4|\h_{n}(z_0)|^2+8|\h_{n+1}(z_0)\h_n(z_0)|} {|\vp_n(\lambda)|^2|z_0-\lambda|^2} \\
& &\leq  \ds \frac {16|\h_{n}(z_0)|^2} {|\vp_n(\lambda)|^2|z_0-\lambda|^2}  \, .
\label{ineq}
\end{array}
\end{equation}

\emph{Remark:} Note that $\h_{n+1}(z_0)-\h_n(z_0)=(1-\overline{\lambda}z_0)\overline{\vp_n(\lambda)}\vp_n(z_0)$, so it is impossible that both $\h_{n+1}(z_0)$ and $\h_n(z_0)$ are zero because $\vp$ has zeros inside the unit circle. \\ 

On the other hand, we observe that
\begin{equation}
\|K_{n-1}(z_0,\cdot)\|=\left( \int_{\T}K_{n-1}(z_0,y)\overline{K_{n-1}(z_0,y)}d\mu(y) \right)^{1/2}=K_{n-1}(z_0,z_0)^{1/2} \, .
\end{equation}
Hence
\begin{equation} \|(z_0- \cdot)K_{n-1}(z_0, \cdot)\|^2
\geq \dist(z_0, \supp(d\mu))^2 K_{n-1}(z_0,z_0) \, .
\smallskip\end{equation}
As a result,
\begin{equation}
\dist(z_0, \supp(d\mu))^2 K_{n-1}(z_0,z_0) \leq \ds \frac {16|\h_{n}(z_0)|^2} {|\vp_n(\lambda)|^2|z_0-\lambda|^2} \, .
\end{equation}
This proves the case when $|\h_{n+1}(z_0)|\leq|\h_n(z_0)|$.

Now suppose $|\h_{n+1}(z_0)|\leq|\h_n(z_0)|$. The proof could be carried out in a similar manner, only that after (\ref{ineq}) all appearances of $\h_{n}$ will be replaced by $\h_{n+1}$.

\end{proof}\medskip

\begin{lemma}\label{dslemma} Suppose $\tau$ is a zero of $\h_n$ which is distinct from $\lambda$. Let $T=\dist(\tau,\supp(d\mu))$, then
\begin{equation}
|z_0-\tau|\geq \ds \frac{|\h_n(z_0)|}{K_{n-1}(z_0,z_0)^{1/2}\|\h_n\|} \quad T
\, .
\end{equation} 
\label{lemma2}
\end{lemma}

\begin{proof}
Since $\tau$ is a zero of $\h_n$, $g(z)=\frac{\h_n(z)}{(z-\tau)}$ is a polynomial of degree $n-1$, so we can express it as
\begin{equation}
\ds \frac {\h_n(z)}{(z-\tau)}=\ds \int_{\T}K_{n-1}(z,y)g(y)d\mu(y) \, .
\end{equation}
By the Schwarz inequality,
\begin{equation}\left| \ds \frac {\h_n(z_0)}{(z_0-\tau)}\right| \leq \|K_{n-1}(z_0,\cdot)\|\|g\| = K_{n-1}(z_0,z_0)^{1/2}\|g\|
\, .
\end{equation}

Also note that $\|g\|=\left\| \frac{\h_n(z)}{(z-\tau)}\right\| \leq \frac{\|\h_n\|}{T}$. \quad 
Therefore,
\begin{equation}
|z_0-\tau|\geq \ds \frac{|\h_n(z_0)|}{K_{n-1}(z_0,z_0)^{1/2}\|\h_n\|} \quad T \, .
\end{equation}
\end{proof}

\begin{proof}[Proof of Theorem \ref{theorem1}] Notice that either one of the following must be true:

\begin{align}
\label{one}|\h_{n+1}(z_0)| & \leq |\h_n(z_0)| \\
|\h_n(z_0)| & \leq |\h_{n+1}(z_0)| \label{two} \, .
\end{align}

We observe that 
\begin{equation}
\|\h_n\| = \|\overline{\vp^*_n(\lambda)} \vp^*_n(y)-\overline{\vp_{n}(\lambda)}\vp_n(y)\|_{L^2(d\mu(y))} \leq 2|\vp_n(\lambda)| \, .
\label{hnnorm}
\end{equation}
If (\ref{one}) is true, combining this with Lemma \ref{lemma1} and Lemma \ref{dslemma}, we obtain that:
\begin{equation}
|z_0-\tau| \geq  \left(\ds \frac{\delta^2 |\vp_n(\lambda)|}{4}\ds \frac{1}{2|\vp_n(\lambda)|} \right) T  = \ds \frac{\delta^2 T}{8} \, .
\end{equation}
Finally, by the triangle inequality,
\begin{equation}
T=\dist(\tau,\supp(d\mu)) \geq \dist(z_0,\supp(d\mu))-|z_0-\tau|=\delta-|z_0-\tau| \, .
\end{equation}
This gives
\begin{equation}|z_0-\tau|\geq \ds \frac{ \delta^2(\delta-|z_0-\tau|) }{8}
\end{equation}
and the result follows.

On the other hand, if (\ref{two}) is true, then instead of (\ref{hnnorm}) we use the definition of $\h_{n+1}$ in (\ref{hndef2}) which will give the same bound of $\|\h_{n+1}\|$ as in (\ref{hnnorm}). Hence the same argument applies to $\h_{n+1}$.

\smallskip

Now consider the special case where $L=\dist(\lambda,\supp(d\mu))>0$. Without loss of generality, suppose (\ref{one}) is true. Since $\tau$ and $\lambda$ are distinct zeros of $\h_n$, we could apply a similar argument as in Lemma {\ref{lemma2}} to $\frac {\h_n(z)}{(z-\tau)(z-\lambda)}$ and obtain the following
\begin{equation} \begin{array} {lll}
|z_0-\tau||z_0-\lambda| & \geq & \ds \frac{|\h_n(z_0)|}{K_{n-2}(z_0,z_0)^{1/2}\|\h_n\|} T L \, . \\
\end{array}
\label{cor}
\end{equation}
Since $K_{n-2}(z_0,z_0)^{1/2} \leq K_{n-1}(z_0,z_0)^{1/2}$, the desired inequality follows.
Now we combine (\ref{cor}) with Lemma $\ref{lemma1}$. The $|z_0-\lambda|$ term cancels on both sides and it gives us
\begin{equation}
|z_0-\tau| \geq \ds \frac{\delta L T}{8} \, .
\end{equation}
Again, we use the triangle inequality on T and the result follows. Clearly, if ({\ref{two}}) is true, we could still apply the same argument to $h_{n+1}$.
\end{proof}\medskip

\section{Proof of Theorem \ref{theorem2}}
\begin{proof} According to the definitions of $\vp_n^*$ and $\p_n^*$,
\begin{eqnarray}
\s_n(z)= \overline{\lambda^n}z^n \vp_n(\lambda)\overline{\p_n(z)}+ \overline{\vp_n(\lambda)} \p_n(z)\\
\h_n(z)= \overline{\lambda^n}z^n \vp_n(\lambda)\overline{\vp_n(z)}- \ol{\vp_n(\lambda)}\vp_n(z)\, .
\end{eqnarray}

If we define for $z \in \T$
\begin{eqnarray}
\sigma_n(z):=\ds \frac{\s_n(z)}{(\overline{\lambda}z)^{n/2}} \label{sigman}\\
\eta_n(z):=\ds \frac{h_n(z)}{i(\overline{\lambda}z)^{n/2}} \label{etan}
\end{eqnarray} with $\rm{Arg}((\overline{\lambda}z)^{1/2}) \in [0,\pi)$, then $\sigma_n$ and $\eta_n$ are real-valued $C^{\infty}$ functions and they have the same zeros as $\s_n$ and $\h_n$ respectively.

To prove the interlacing condition of Theorem \ref{theorem2}, it suffices to prove the following:
\begin{equation}
\ds \frac{d \eta_n(e^{i\theta})}{d\theta} \sigma_n(e^{i\theta})<0 \text{ at every zero } e^{i \theta} \text{ of } \eta_n(z) \, .
\label{equivstatement}
\end{equation}

We shall prove condition (\ref{equivstatement}) for $n+1$.

Suppose $\zeta$ is a zero of $\h_{n+1}$. By (\ref{formulah}), $\h_{n+1}$ could be expressed by the reproducing kernel. Hence $\eta_{n+1}$ can be represented as \begin{equation}
\eta_{n+1}(z)=\frac{1}{i(\ol{\lambda}z)^{(n+1)/2}}\ds \frac{-\overline{\lambda \vp_n(\lambda)}}{\overline{\vp_n(\zeta)}}(z-\zeta)\ds \sum_{j=0}^{n}\vp_j(z)\overline{\vp_j(\zeta)}
\, .
\label{cmp}
\end{equation}
The constant $\frac{-\overline{\lambda \vp_n(\lambda)}}{\overline{\vp_n(\zeta)}}$ is obtained by comparing the leading coefficients of the right hand side of (\ref{cmp}) and that of $\h_{n+1}$ when expressed in terms of (\ref{hndef2}).

As a result, the derivative of $\eta_{n+1}$ at $\zeta$ is
\begin{equation}\begin{array}{ll}
\ds \frac{d \eta_{n+1}}{dz}(\zeta)& = \ds \lim_{z\rightarrow \zeta} \ds \frac{\eta_{n+1}(z)-\eta_{n+1}(\zeta)}{z-\zeta} \\ & = \ds \lim_{z\rightarrow \zeta} \ds \frac{\eta_{n+1}(z)}{z-\zeta} \\ & = \ds \frac{-\overline{\lambda \vp_n(\lambda)}}{i \overline{\vp_n(\zeta)}} \left( \frac{\lambda}{\zeta}\right)^{\frac{n+1}{2}}K_n(\zeta,\zeta) \, .
\end{array}
\end{equation}

Let $\zeta=e^{i \theta}$ and $z=e^{i\omega}$. By the chain rule,
\begin{equation} \begin{array}{ll}
\ds \frac{d \eta_{n+1}}{d \omega}(\theta) & = i \zeta \ds \frac{d \eta_{n+1}}{dz}(\zeta) \\
& = \ds -\frac{\overline{\vp_n(\lambda)}}{\overline{\vp_n(\zeta)}} \left( \frac{\lambda}{\zeta}\right)^{\frac{n-1}{2}}K_n(\zeta,\zeta) \, .
\end{array}
\end{equation}

Now we go back to $\ds \frac{d \eta_n(e^{i\theta})}{d\theta} \sigma_n(e^{i\theta})$ and compute:
\begin{equation}\begin{array}{ll}
& \ds \frac{d \eta_{n+1}(e^{i\theta})}{d\theta} \sigma_{n+1}(e^{i\theta}) \\
= & \ds -\frac{\overline{\vp_n(\lambda)}}{\overline{\vp_n(\zeta)}} \left( \frac{\lambda}{\zeta}\right)^{n} K_n(\zeta,\zeta)
\left( \ol{\vp_{n}^*(\lambda)}\p_{n}^*(\zeta)+\ol{\lambda}\zeta \ol{\vp_{n}(\lambda)}\p_{n}(\zeta) \right) \\
 = & \ds - \left( \frac{\lambda}{\zeta}\right)^{n} K_n(\zeta,\zeta) \left(
|\vp_n(\lambda)|^2 \left( \frac{\zeta}{\lambda}\right)^{n} \ds \frac{\ol{\p_n(\zeta)}}{\ol{\vp_n(\zeta)}}
+\ol{\lambda}\zeta \ds \frac{\overline{\vp_n(\lambda)}}{\overline{\vp_n(\zeta)}} \ol{\vp_n(\lambda)}\p_n(\zeta)
\right) \, .
\end{array}
\label{eqnabv}
\end{equation}

Recall that $\eta_{n+1}(\zeta)=0$, which implies that
\begin{equation}
\ds \frac{\overline{\vp_n(\lambda)}}{\overline{\vp_n(\zeta)}} = \ds \frac{\vp_n(\lambda)}{\vp_n(\zeta)}\left( \frac{\zeta}{\lambda}\right)^{n-1} \, .
\end{equation}

We then apply this onto the second part of the summand in (\ref{eqnabv}):
\begin{align}
(\ref{eqnabv}) & =  
\ds - \left( \frac{\lambda}{\zeta}\right)^{n} K_n(\zeta,\zeta) \left(
|\vp_n(\lambda)|^2 \left( \frac{\zeta}{\lambda}\right)^{n} \ds \frac{\ol{\p_n(\zeta)}}{\ol{\vp_n(\zeta)}} + 
\left( \ds \frac{\zeta}{\lambda} \right)^n \ds \frac{\vp_n(\lambda)}{\vp_n(\zeta)} \ol{\vp_n(\lambda)}\p_n(\zeta)\right) \notag \\
&=  \ds - K_n(\zeta,\zeta) |\vp_n(\lambda)|^2 \left( \ds \frac{\ol{\p_n(\zeta)}}{\ol{\vp_n(\zeta)}}+ \frac{\p_n(\zeta)}{\vp_n(\zeta)}\right) \notag \\ 
& = \ds - K_n(\zeta,\zeta) \ds \left| \frac{\vp_n(\lambda)}{\vp_n(\zeta)} \right|^2
\left( \ol{\p_n(\zeta)}\vp_n(\zeta)+\ol{\vp_n(\zeta)}\p_n(\zeta) \right)
\label{laeqn} \, .
\end{align}

Now we use a formula that relates $\vp_n$ and $\p_n$ (see Chapter 3.2 in \cite{simon}):
\begin{equation}
\ol{\p_n(z)}\vp_n(z)+\ol{\vp_n(z)}\p_n(z) = 2 \text{ in }\T
\label{relformula} \, .
\end{equation}

We apply (\ref{relformula}) to (\ref{laeqn}). This gives us the result that at any zero $\zeta$ of $\eta_{n+1}$:
\begin{equation}
\ds \frac{d \eta_{n+1}(e^{i\theta})}{d\theta} \sigma_{n+1}(e^{i\theta}) = (\ref{eqnabv})=\ds - 2 K_n(\zeta,\zeta) \ds \left| \frac{\vp_n(\lambda)}{\vp_n(\zeta)} \right|^2 <0 \, .
\end{equation}

The interlacing theorem is proven.
\end{proof}
\smallskip

\section{Proof of lemma \ref{mainlemma}}
We prove Lemma \ref{mainlemma} by stating several lemmas which are similar to those in the proof of Theorem \ref{theorem1}.
\begin{lemma}
\label{lemma5} Suppose $\dt=\dist(z_0,\supp(d\nu))>0$ and $\kt_n(x,y)=\ds \sum_{j=0}^{n}\p_j(x)\ol{\p_j(y)}$ is the reproducing kernel with respect to the measure $\nu$. Then
\begin{equation}  \ds \left|\frac{\s_i(z_0)}{\kt_{n-1}(z_0,z_0)^{1/2}} \right| \geq \frac{1}{4}|\vp_n(\lambda)| |z_0-\lambda|\dt 
\end{equation} where $i= \left\{ \begin{array}{lll} n & \text{ if }\hspace{1cm} |\s_{n+1}(z_0)| & \leq |\s_n(z_0)| \\ n+1 & \text{ if }\hspace{1cm} |\s_{n}(z_0)| & \leq |\s_{n+1}(z_0)| \end{array} \right. $ .
\end{lemma}

\begin{proof} The proof is essentially the same as the one of Lemma \ref{lemma1}, except for a few differences. The $L^2$ norm here refers to the one taken with respect to $\nu$ and $\h_n$ is replaced by $\s_n$.

It is also worth noting that by the definition of $\s_n$ in (\ref{sndef3}),
\begin{equation}
\s_{n+1}(z)-\s_n(z) = -(1-\overline{\lambda}z)\overline{\vp_{n}(\lambda)}\p_{n}(z) \not = 0 \text{ on } \T \, .
\end{equation}
As a result,
\begin{equation}|\p_{n}(z_0)| = \ds \left| \frac{\s_{n+1}(z_0)-\s_n(z_0)}{(z_0-\lambda)\vp_{n}(\lambda)} \right| 
\end{equation}
which allows us to proceed in the same way as in the proof of Lemma \ref{lemma1}.
\end{proof}
\smallskip
\begin{lemma} \label{lemma6} Suppose $\ttau$ is a zero of $\s_n$. Let $\tT=\dist(\ttau,\supp(d\nu))$, then
\begin{equation}
|z_0-\ttau|\geq \ds \frac{|\s_n(z_0)|}{\kt_{n-1}(z_0,z_0)^{1/2}\|\s_n\|_{L^2(d\nu)}} \tT  \, . \end{equation}
\end{lemma}
The proof of this lemma is omitted because it resembles that of Lemma \ref{lemma2}.
\medskip

Finally, we state the following lemma relating the support of $\mu$ and $\nu$:
\begin{lemma} \label{lastlemma} Suppose $z_0$ is an isolated point in the support of $\mu$. Then
\begin{equation}
\dt=\dist(z_0,\supp(d\nu))>0 \, .
\end{equation}
\end{lemma}
The reader could refer to Chapter 3.2, p225 of \cite{simon} for the proof.

\medskip

Next, we are going to finish the proof of Lemma \ref{mainlemma}.
\begin{proof} Suppose $z_0$ is an isolated point in the support of $d\mu$ which is distinct from $\lambda$. By Lemma \ref{lastlemma}, $\dist(z_0,\supp(d\nu))>0$.

Either $|\s_n(z_0)|\geq |\s_{n+1}(z_0)|$ or $|\s_n(z_0)|\leq |\s_{n+1}(z_0)|$ is true.
Without loss of generality, we assume that $|\s_n(z_0)| \geq |\s_{n+1}(z_0)|$ and use Lemma \ref{lemma5}.

Furthermore, we observe that
\begin{equation}
\|\s_n\| \leq 2|\vp_n(\lambda)|\|\p_n\|_{L^2(d\nu)}=2|\vp_n(\lambda)|
\end{equation}
Then we combine these results to get
\begin{equation}
|z_0-\ttau| \geq \ds \frac{|z_0-\lambda| \dt \tT}{8} \, .
\end{equation}

Finally, we apply the triangle inequality to $\tT$:
\begin{equation}
\tT=\dist(\ttau,\supp(d\nu)) \geq \dist(z_0,\supp(d\nu))-|z_0-\ttau|=\dt-|z_0-\ttau| \, .
\end{equation}
This gives us the following inequality which finishes the proof:
\begin{equation}
|z_0-\ttau| \geq \ds \frac{\dt^2 |z_0-\lambda|}{8+|z_0-\lambda|\dt} \, .
\label{bound}
\end{equation}
\end{proof}

\section{Proof of Theorem \ref{theorem3}}
\begin{proof} By Lemma \ref{mainlemma}, inside the ball $B(z_0,\rt)$ either $\s_n$ or $\s_{n+1}$ (or both) has no zero inside, with $\rt$ given by (\ref{bound}) above. Without loss of generality, we assume that $\s_n$ does not have zeros inside. By Theorem \ref{theorem2} the zeros of $\h_n$ and $\s_n$ interlace, therefore $\h_n$ cannot have more than two zeros inside $B(z_0,\rt)$.
\end{proof}
\smallskip

\section{Acknowledgements}
I would like to thank Professor Barry Simon for his suggesting this problem, as well as his time for many very helpful discussions and email communications. I would also like to thank Cherie Galvez for her editorial advice as well as her help with LaTeX.

\end{document}